\documentclass[reqno,11pt]{amsart}
\usepackage{ifthen, citesort}

\usepackage{latexsym}
\usepackage{graphicx}
 \usepackage{amsmath}
 \usepackage{amsfonts}
 \usepackage{amssymb}
 \usepackage{amstext}
 \usepackage{amsbsy}
 \usepackage{amsthm}
 \usepackage{amscd}
\usepackage{amsfonts}
\usepackage{amssymb}
\usepackage{amsmath}
\usepackage{latexsym, amsfonts}
\usepackage{portland}
\usepackage{longtable}
\usepackage[T1]{fontenc}

\newcommand{\psl}{PSL\ensuremath{(2,\mathbb{C})\ }}

\newcommand{\nospacepsl}{PSL\ensuremath{(2,\mathbb{C})}}
\newcommand{\nospaceslc}{SL\ensuremath{(2,\mathbb{C})}}
\newcommand{\tr}{\mathrm{tr}}
\newcommand{\QQ}{\ensuremath{\mathbb{Q}}}
\newcommand{\ZZ}{\ensuremath{\mathbb{Z}}}
\newcommand{\h}{\ensuremath{\mathbb{H}^3} }
\newcommand{\pf}{{\bf Proof:  }}
\newcommand{\done}{\hfill $\square$ \newline}
\newcommand{\Vol}{\rm{Vol}}
\newcommand{\PSL}{\rm{PSL}}

\newtheorem{thm}{Theorem}
\newtheorem*{thmar}{Theorem}

\newtheorem{prop}{Proposition}

\newtheorem{rmk}{Remark}

\newtheorem{conj}{Conjecture}

\begin{document}

%%%%%%%%%TITLES ECT%%%%%%%%%%%%%%%%%%%%%%%%%%%%%%%%%%%%
\title{Exceptional regions and associated exceptional hyperbolic 3-manifolds}

\author{Abhijit Champanerkar, Jacob Lewis, Max Lipyanskiy and Scott Meltzer
\\ With an appendix by Alan W. Reid}
\address{Department of Mathematics and Statistics, University of South Alabama, Mobile, AL 36688}
\email{achampanerkar@jaguar1.usouthal.edu}

%\author{Jacob Lewis}

\address{Department of Mathematics,
University of Washington,
Seattle, WA 98195-4350}
\email{jacobml@math.washington.edu}

%    Information for third author
%\author{ Max Lipyanskiy}

\address{Department of Mathematics,  Massachusetts Institute of Technology
    77 Massachusetts Avenue
    Cambridge, MA 02139-4307}
\email{mlipyan@math.mit.edu}

%    Information for fourth author
%\author {Scott Meltzer}

\address{4515 Woodland Park Ave North, Seattle, WA 98103}
\email{smeltzer2003@yahoo.com}

% Information for author of appendix
%\author{\\ With an appendix by Alan Reid}

\address{Department of Mathematics, University of Texas, Austin, TX 78712}
\email{areid@math.utexas.edu}

%\leftheadtext\nofrills{

%    General info
\renewcommand{\subjclassname}{%
  \textup{2000} Mathematics Subject Classification}
\subjclass[2000]{Primary 57M50; Secondary 57N10}

\date{}

%\dedicatory{}

%\keywords{}

\begin{abstract}

A closed hyperbolic 3-manifold is exceptional if its shortest geodesic
does not have an embedded tube of radius $\ln(3)/2$. D. Gabai,
R. Meyerhoff and N. Thurston identified seven families of exceptional
manifolds in their proof of the homotopy rigidity theorem. They identified
the hyperbolic manifold known as Vol3 in the literature as the exceptional
manifold associated to one of the families.  It is
conjectured that there are exactly 6 exceptional manifolds.
 We find hyperbolic 3-manifolds, some from the SnapPea's
census of closed hyperbolic 3-manifolds, associated to 5 other
families. Along with the hyperbolic 3-manifold found by Lipyanskiy
associated to the seventh family we show that any exceptional manifold
is covered by one of these manifolds. We also find their group coefficient
fields and invariant trace fields.

\end{abstract}

\maketitle

\markboth{A. Champanerkar, J. Lewis, M. Lipyanskiy, S. Meltzer and A. W. Reid}
{Exceptional regions and associated exceptional hyperbolic 3-manifolds }

%%%%%%%%%%%%%%%%%%%%%%%%%%%%%%%%%%%%%%%%%%%%%%%%%%%%%%%%%%

%%%%%%%%%%%%%PAPER%%%%%%%%%%%%%%%%%%%%%%%%%%%%%%%%

\section{Introduction}

A closed hyperbolic 3-manifold is \textit{exceptional} if its shortest
geodesic does not have an embedded tube of radius
$\ln(3)/2$. Exceptional manifolds arise in the proof of the rigidity
theorem proved by D. Gabai, R. Meyerhoff and N. Thurston in
\cite{gmt}. It is conjectured that there are exactly 6 exceptional
manifolds.

%% It is shown in \cite{gmt} and \cite{jr01} that Vol3, the
%% closed hyperbolic 3-manifold with conjecturally the third smallest
%% volume, is exceptional. We show that the manifolds $v2678(2,1)$,
%% $s778(-3,1)$ and $s479(-3,1)$ from SnapPea's census of closed
%% hyperbolic 3-manifolds \cite{snappea} and the manifold described in
%% Section 5 which we call $N_4$ are exceptional. Lipyanskiy has
%% described a sixth exceptional manifold in \cite{lipyan}. We show:

%% \begin{thm}
%% Let $N$ be an exceptional manifold. Then $N$ is covered by one of the
%% above 6 exceptional manifolds.
%% \end{thm}

Let $N$ be a closed hyperbolic 3-manifold and $\delta$ be the shortest
geodesic in $N$. If $\delta$ does not have an embedded tube of radius
$\ln(3)/2$ then there is a two-generator subgroup $G$ of $\pi_1(N)$
such that $\h/G$ also has this property. Assume $G$ is generated by
$f$ and $w$, where $f \in \pi_1(N)$ is a primitive hyperbolic isometry
whose fixed axis $\delta_0 \in \mathbb{H}^3$ projects to $\delta$ and
$w\in \pi_1(N)$ is a hyperbolic isometry which takes $\delta_0$ to its
nearest translate. Thus, it is necessary to study two-generator
subgroups $\Gamma$ of PSL$(2,\mathbb{C})$ with the property that one
of the generators is the shortest geodesic in $\h/\Gamma$ and the
distance from its nearest translate is less than $\ln(3)$.

The space of two-generator subgroups of \psl is analysed in the proof
of the rigidity theorem of \cite{gmt}. The rigidity theorem is proved
using Gabai's theorem \cite{gabai}, which states that the rigidity
theorem is true if some closed geodesic has an embedded tube of radius
$\ln(3)/2$. The authors of \cite{gmt} show that this holds for all but
seven  exceptional families of closed hyperbolic 3-manifolds.  These
seven families are handled separately. The seven families are obtained
by parametrizing the space of two-generator subgroups of \psl by a
subset of $\mathbb{C}^{3}$, dividing the parameter space into about a
billion regions and eliminating all but seven regions. These seven
regions correspond to the seven exceptional families and are known as
\textit{exceptional regions}. They are denoted by $X_i$ for $i=0,
\ldots, 6$ and described as boxes in $\mathbb{C}^{3}$. For example,
for the region $X_3$ see Table \ref{rangeforx3}.

\begin{table}[!h]
\begin{center}
\begin{tabular}{|l|c|r|}

\hline Parameter & Range $Re(Parameter)$ & Range $Im(Parameter)$\\
\hline $L'$ & 0.58117 to 0.58160 & -3.31221 to -3.31190\\ \hline $D'$
& 1.15644 to 1.15683 &-2.75628 to -2.75573\\ \hline $R'$ & 1.40420 to
1.40454 & -1.17968 to -1.17919\\ \hline
\end{tabular}\\
\end{center}
\vspace{.1in}

\caption{Parameter ranges for the region $X_3$}
\label{rangeforx3}
\end{table}

% These seven families are described by regions in $\mathbb{C}^{3}$ using a parametrization of the
%space of two-generator subgroups by a subset of $\mathbb{C}^{3}$ described in \cite{gmt}.
%These seven regions are known as \textit{exceptional regions} and denoted by $X_i$ for
%$i=0, \ldots, 6$. For example, for the region $X_3$ see Table \ref{rangeforx3}.

 A \textit{quasi-relator} in a region is a word in $f$, $w$, $f^{-1}$
and $w^{-1}$ that is close to the identity throughout the region and
experimentally  appears to converge to the identity at some
point. Table \ref{relators} gives the two quasi-relators specified for
each region $X_i$ in \cite{gmt}. The group $G_i= \langle f,\ w|
r_1(X_i),\ r_2(X_i)\rangle$, where $ r_1(X_i),\ r_2(X_i)$ are the
quasi-relators for $X_i$, is called the \textit{marked group} for the
region $X_i$. It follows from \cite{gmt} that any exceptional manifold
has a two-generator subgroup of its fundamental group whose parameter
lies in one of the exceptional regions.

\begin{table}[!h]
\begin{center}
\begin{tabular}{|l|p{4in}|}
\hline Region & Quasi-Relators  \\
\hline
$X_0$ &
$r_1= fwf^{-1}w^2f^{-1}wfw^2$\\
& $r_2=f^{-1}wfwfw^{-1}fwfw$\\

\hline
$X_1$ &

$r_1=f^{-2}wf^{-1}w^{-1}f^{-1}w^{-1}fw^{-1}f^{-1}w^{-1}f^{-1}wf^{-2}w^2$\\
 & $r_2=f^{-2}w^2f^{-1}wfwfw^{-1}fwfwf^{-1}w^2$ \\

\hline $X_2$ &

$r_1=f^{-1}wfwfw^{-1}f^2w^{-1}fwfwf^{-1}w^2$ \\
& $r_2=f^{-2}wf^{-2}w^2f^{-1}wfwfwf^{-1}w^2$ \\

\hline

$X_3$ &

$r_1=f^{-2}wfwf^{-2}w^2f^{-1}w^{-1}f^{-1}wf^{-1}w^{-1}(fw^{-1}f^{-1}w^{-1}f)^2$
\\ &
$w^{-1}f^{-1}wf^{-1}w^{-1}f^{-1}w^2$
\\ &
$r_2=f^{-2}wfwf^{-1}wf(w^{-1}fwfw^{-1})^2fwf^{-1}wfwf^{-2}w^2$
\\ &
$f^{-1}w^{-1}f^{-1}w^2$\\

\hline
$X_4$ &
$r_1= f^{-2}wfwf^{-1}(wfw^{-1}f)^2wf^{-1}wfwf^{-2}w^2$\\ & $(f^{-1}w^{-1}f^{-1}w)^2w $
\\
& $r_2=f^{-1}(f^{-1}wfw)^2f^{-2}w^2f^{-1}w^{-1}f^{-1}w(f^{-1}w^{-1}fw^{-1})^2$
\\ &$f^{-1}wf^{-1}w^{-1}f^{-1}w^2 $\\
\hline
$X_5$ &

$r_1=f^{-1}wf^{-1}w^{-1}f^{-1}wf^{-1}wfwfw^{-1}fwfw$\\
&  $r_2=f^{-1}wfwfw^{-1}fw^{-1}f^{-1}w^{-1}fw^{-1}fwfw$  \\
\hline

$X_6$ &

$r_1=f^{-1}w^{-1}f^{-1}wf^{-1}w^{-1}f^{-1}w^{-1}fw^{-1}fwfw^{-1}fw^{-1}$\\
&  $r_2=f^{-1}w^{-1}fw^{-1}fwfwf^{-1}wfwfw^{-1}fw^{-1}$  \\

\hline

\end{tabular}

\end{center}
\vspace{.1in}
\caption{Quasi-relators for all the regions}
\label{relators}
\end{table}

Let $\rho(x,y)$ denote the hyperbolic distance between $x,\ y \in \h$.
 For an isometry $f$ of $\h$ define
$Relength(f)={\rm inf} \{ \rho(x,f(x))|\  x \in \h \}$.
 Let $T$ consists of those parameters corresponding to the groups
 $\{G,\ f,\ w\}$ such that $Relength(f)$ is the shortest element of
 $G$ and the distance between the axis of $f$ and its nearest
 translate is less than $\ln(3)$. Let $S=exp(T)$. In \cite{gmt} the authors
 made the following conjecture.

\begin{conj} \label{gmtconj}
Each sub-box $X_i$, $0\leq i \leq 6$ contains a unique element $s_i$ of S.
 Further if $\{G_i,f_i,w_i\}$ is the marked group associated to $s_i$ then $N_i=\h/G_i$
is a closed hyperbolic 3-manifold with the following properties
\begin{enumerate}
\item[(i)] $N_i$ has fundamental group $\langle f,\ w | r_1(X_i),\
r_2(X_i) \rangle $ where $r_1$ and $r_2$ are the quasi-relators
associated to the box $X_i$.
\item[(ii)]$N_i$ has a Heegaard genus 2 splitting realizing the above
group presentation.
\item[(iii)]$N_i$ nontrivially covers no manifold.
\item[(iv)]$N_6$ is isometric to $N_5$.
\item[(v)]If $(L_i,D_i,R_i)$ is the parameter in $T$ corresponding to
$s_i$, then $L_i,\ D_i,$ $R_i$ are related as follows: \\ For $X_0, \
X_5,\ X_6$: $L=D$, $R=0$.\\ For $X_1,\ X_2,\ X_3,\ X_4$: $R=L/2$.
\end{enumerate}

\end{conj}

D. Gabai, R. Meyerhoff and N. Thurston \cite{gmt} proved that Vol3,
the closed hyperbolic 3-manifold with conjecturally the third smallest
volume, is the unique exceptional manifold associated to the
region $X_0$. Jones and Reid \cite{jr01} proved that Vol3 does not
nontrivially cover any manifold and that the exceptional manifolds
associated to the regions $X_5$ and $X_6$ are isometric.

In this paper we investigate the seven exceptional regions and the
associated exceptional hyperbolic 3-manifolds.
 In Section 3, using Newton's method  for finding roots
of polynomials in several variables, we solve the  equations obtained
from the entries of the quasi-relators to very high precision. Then,
using the program  PARI-GP, \cite{parigp} we find entries of the
generating matrices as  algebraic numbers, find the group  coefficient
fields and verify with exact arithmetic that the quasi-relators are
relations for all the regions. We also find the invariant trace fields
for all the groups verifying and extending the data given in
\cite{jr01}.  In Section 4  we show that the manifolds
$v2678(2,1)$, $s778(-3,1)$ and $s479(-3,1)$
from SnapPea's census of closed hyperbolic 3-manifolds \cite{snappea}
are the exceptional manifolds associated to the regions $X_1,\ X_2,\
X_5$ and $X_6$ respectively.  We also show that their fundamental
groups are isomorphic to the marked groups $G_i$ for $i=1,2,5,6$.
In Section 5  we find an exceptional manifold associated to the region
$X_4$ and show that its fundamental group is isomorphic to the marked
group $G_4$. This manifold, which we denote by $N_4$, is commensurable
to the SnapPea census manifold $m369(-1,3)$.  Lipyanskiy has
 described a sixth exceptional manifold in \cite{lipyan}.

In Section 6, using
Groebner bases we show that the quasi-relators have a unique solution
in every region. Let $N_0=$Vol3, $N_1=v2678(2,1),\ N_2=s778(-3,1)$,
$N_3$ be the exceptional manifold associated to $X_3$ found in \cite{lipyan},
$N_4$ be exceptional manifold associated to $X_4$ found in Section 5 and
$N_5=s479(-3,1)$. We shall prove the following theorem.

\begin{thm}Let $N$ be an exceptional manifold. Then $N$ is covered by
$N_i$ for  some $i=0,1,2,3,4,5$.
\end{thm}

\noindent
\textbf{Acknowledgments:} We thank the Columbia
VIGRE and I. I. Rabi programs for supporting the summer undergraduate
research projects. We thank Walter Neumann for his continuous
guidance, encouragement and support. We thank Alan Reid for his valuable
comments and suggestions.

 \section{Invariant trace fields and 2-generator subgroups }

  Two finite volume, orientable, hyperbolic 3-manifolds are said to be
\textit{commensurable} if they have a common finite-sheeted cover.
Subgroups $G, G'\subset {\rm PSL}(2,\mathbb{C})$ are \textit{commensurable} if there exists
$g \in {\rm PSL}(2,\mathbb{C})$ such that $g^{-1}Gg \cap G'$ is a finite index
subgroup of both $g^{-1}Gg$ and $G'$. It follows by Mostow rigidity
that finite volume, orientable, hyperbolic 3-manifolds are commensurable if
and only if their fundamental groups are commensurable as subgroups of
\nospacepsl. Let $G$ be a group of covering transformations and let $\tilde{G}$
be its preimage in \nospaceslc. It is shown in \cite{macbeath83} that the
traces of elements of $\tilde{G}$
generate a number field $\mathbb{Q}(\tr G)$ called the \textit{trace field}
of G. The \textit{invariant trace field} $k(G)$ of $G$ is defined as the
intersection of all the fields $\mathbb{Q}(\tr H)$ where $H$ ranges over
all finite index subgroups of $G$. The definition already makes clear that
the $k(G)$ is a commensurability invariant. In \cite{reid90} Alan Reid proved
the following result.
\begin{thm} The invariant trace field $k(G)$ is equal to
$$\mathbb{Q}(\{\tr^2(g):g \in G\})=\mathbb{Q}(\tr G^{(2)}), $$
where $G^{(2)}$ is the finite index subgroup of $G$ generated by
squares $\{g^2:g\in G\}$.
\end{thm}

From Corollary 3.2 of \cite{hlm92} the invariant trace field of
a 2-generator group $\langle f, w | r_1, r_2 \rangle $ is generated by
$\tr(f^2), \tr(w^2)$ and $\tr(f^2w^2)$.

Using trace relations (see Theorem 4.2 of \cite{cghn00}) and Corollary 3.2
of \cite{hlm92}, which says that $\QQ(\tr G^{(2)})=\mathbb{Q}(\tr G^{SQ})$
where $G^{SQ}=\langle g_1^2, \dots g_n^2 \rangle$ with $g_i$'s generators
of $G$ such that $\tr(g_i)\neq 0$, we see that the invariant trace field of
a 2-generator group $\langle f, w | r_1, r_2 \rangle $ is generated by
$\tr(f^2), \tr(w^2)$ and $\tr(f^2w^2)$.

As described in the Introduction, a marked group $G$ is generated by $f$ and
$w$, where $f$ and $w$ are seen as covering transformations such that $f$
represents the shortest geodesic and $w$ takes the axis of $f$ to its
nearest translate. Let \h denote the upper half space model of hyperbolic 3-space
and the sphere at infinity be the $x-y$ plane. Conjugate $G$ so that the axis of $f$ is the geodesic
line $B_{(0,\infty)}$ in \h with end points $0$ and $\infty$ on the sphere at infinity and the
geodesic line perpendicular to $w^{-1}(B_{(0,\infty)})$
and $B_{(0,\infty)}$ (orthocurve) lies on the geodesic line $B_{(-1,1)}$ in \h with endpoints
$-1$ and $1$ on the sphere at infinity. We can parametrize such a marked group with $3$ complex
numbers $L,D$ and $R$ where $f$ is an $L$ translation of
$B_{(0,\infty)}$ and $w$ is a $D$ translation of $B_{(-1,1)}$ followed by
an $R$ translation of $B_{(0,\infty)}$. We can write matrix representatives
for $f$ and $w$ using the exponentials $L', D'$ and $R'$ of $L,D$ and $R$
respectively (see Chapter 1 of \cite{gmt}). We have
\begin{equation}
f=
\begin{pmatrix}
\sqrt{L'}&0 \\
0&1/\sqrt{L'}
\end{pmatrix},
\hspace{0.2in}
w=
\begin{pmatrix}
\sqrt{R'}*ch & \sqrt{R'}*sh \\
sh/ \sqrt{R'}& ch/ \sqrt{R'}
\end{pmatrix}
\label{fandw}
\end{equation}

where $ch=(\sqrt{D'}+1/\sqrt{D'})/2$ and $sh=(\sqrt{D'}-1/\sqrt{D'})/2$
We can write down the generators of the invariant trace field of $G$ in
terms of $L',D'$ and $R'$ as follows
\begin{eqnarray}
{\rm tr}(f^2) &=& L'+\frac{1}{L'}\\
{\rm tr}(w^2)&=&\displaystyle{\frac{1}{4}\Bigl[\Big(R'+\frac{1}{R'}+2\Big)\Big(D'+\frac{1}{D'}+2\Big)-8\Bigr]}\\
\nonumber {\rm tr}(f^{2}w^2) &=& \displaystyle{
\frac{1}{4}\Bigl[\Big(D'+\frac{1}{D'}+2\Big)\Big(R'L'+\frac{1}{R'L'}\Big) } \\
 & &\ \ + \displaystyle{ \Big(D'+\frac{1}{D'}-2\Big)\Big(L'+\frac{1}{L'}\Big)\Bigr]}
\end{eqnarray}

 \section{Guessing the algebraic numbers and exact arithmetic}

 In this section we find the marked groups
for the regions as subgroups of \psl with algebraic entries and find
their invariant trace fields. In \cite{gmt} parameter ranges for the
seven regions are specified. For example, for the parameter range for region
$X_3$ see Table \ref{rangeforx3}.
% \begin{table}
% \begin{center}
% \begin{tabular}{|l|c|r|}

% \hline Parameter & Range $Re(Parameter)$ & Range $Im(Parameter)$\\
% \hline $L'$ & 0.58117 to 0.58160 & -3.31221 to -3.31190\\
% \hline $D'$ & 1.15644 to 1.15683 &-2.75628 to -2.75573\\
% \hline $R'$ & 1.40420 to 1.40454 & -1.17968 to -1.17919\\
% \hline
% \end{tabular}\\
% \end{center}
% \caption{Parameter ranges for the region $X_3$}
% \label{rangeforx3}
% \end{table}

We solve for $L', D'$ and $R'$ such that
the quasi-relators are actually relations in the group. We obtain
eight equations in three complex variables out of which three are
independent and we use Newton's method with the parameter range as
approximate solutions to find high precision solutions, e.g 100
significant digits, for the parameters satisfying the equations.
This allows us to compute $ a = \sqrt{L'} $, $ b = \sqrt{R'} $, $ c = \sqrt{D'} $
and $\tr(f^2),\ \tr(w^2),\ \tr(f^2w^2)$ to high precision. Once
the numbers are obtained to high precision we use the algdep() function
of the PARI-GP package \cite{parigp} to guess a polynomial over the
integers that has the desired number as a root. Although the algdep()
function cannot prove that the guess is in fact correct, we
prove this by using the guessed values to perform exact arithmetic and verify
the relations.

 Once we obtain $ a = \sqrt{L'} $, $ b = \sqrt{R'} $, and $ c = \sqrt{D'} $ as
roots of polynomials we find a primitive element which generates the field that contains all the
three numbers.

% Once $ a = \sqrt{L'} $, $ b = \sqrt{R'} $, and $ c = \sqrt{D'} $ have been
% found to high precision, PARI-GP can be used to guess minimal polynomials
% for $a$, $b$, and $c$. Then, if a field extension can be found containing
% $a$, $b$, and $c$, they can be represented exactly in PARI-GP.

For the regions $ X_0 $, $ X_5 $, and $ X_6 $, $a$, $b$, and $c$ are all contained
in $\mathbb{Q}(a)$. By expressing the matrix entries as algebraic numbers
one can verify the relations directly. For example, for $ X_0 $, the
minimal polynomial for $a$ and $c$ is $ x^8 + 2 x^6 + 6 x^4 + 2 x^2 + 1 $,
and $b=1$, so we can express $a$, $b$, and $c$ as follows:
\begin{verbatim}
a = Mod(x, x^8 + 2*x^6 + 6*x^4 + 2*x^2 + 1),
b = Mod(1, x^8 + 2*x^6 + 6*x^4 + 2*x^2 + 1),
c = Mod(x, x^8 + 2*x^6 + 6*x^4 + 2*x^2 + 1).
\end{verbatim}
Then, using the formulae of Section 2, PARI-GP calculates the quasirelators
exactly as follows:
\begin{verbatim}
[Mod(1, x^8 + 2*x^6 + 6*x^4 + 2*x^2 + 1) 0]
[0 Mod(1, x^8 + 2*x^6 + 6*x^4 + 2*x^2 + 1)].
\end{verbatim}

Thus, exact arithmetic verifies rigorously that the $L'$, $D'$, and $R'$
which were calculated for $ X_0 $ using Newton's method are correct and the
quasi-relators are in fact relators.

In general, the group coefficient field can have arbitrary index over the trace field.
In order to keep the degree of the group coefficient field low we
 follow the method described in \cite{lipyan}. Given that $f$,
$w$ are generic ($fw-wf$ is nonsingular), if $f_2, w_2$ are
any matrices in SL$(2,\mathbb{C})$ such that $\tr(f_2)=\tr(f),\ \tr(w_2)=\tr(w)$
and
$\tr(f_2^{-1}w_2)=\tr(f^{-1}w)$ then the two pairs are conjugate. \\
\indent Let ${\rm tr_1} = \tr(f)$, ${\rm tr_2} = \tr(w)$, ${\rm tr_3} = \tr({f^{-1}}w)$.
  Furthermore let
\begin{equation}
f_2=
\begin{pmatrix}
0&1\\
-1&{\rm tr_1}
\end{pmatrix},
\hspace{0.2in}
w_2=
\begin{pmatrix}
z&0 \\
{\rm tr_1}*z-{\rm tr_3}&{\rm tr_2}-z
\end{pmatrix},
\label{fandwnew}
\end{equation}

% \begin{equation}
% f_2= \left(\begin{array}{cc} 0&1 \\
% -1&{\rm tr_1}
% \end{array}\right)
% \end{equation}

% \begin{equation}
% w_2=\left(\begin{array}{cc} z&0 \\
% {\rm tr_1}*z-tr_3&tr_2-z
% \end{array}\right)
% \end{equation}
\noindent
where $({\rm tr_2}-z)*z=1$. Then the pair $(f_2,w_2)$ is
conjugate to $(f,w)$. The coefficients
of the original $f$ and $w$ may have arbitrary index over the trace
field but in this form the entries of the matrices
are in an at most degree two extension of the trace field.
Table \ref{groupfields} and \ref{groupentries}
 displays the computation of $z$, ${\rm tr_1}$, ${\rm tr_2}$ and ${\rm tr_3}$ for all
regions. In all cases $z$ is the primitive element and ${\rm tr}_i \in
\mathbb{Q}(z)$. One easily verifies the relations using the tables. We have the
following theorem.
\begin{thm}
The marked groups $G_i$ are 2-generator subgroups of ${\rm PSL}(2,\mathbb{C})$ with entries in the
number-fields as given in Tables \ref{groupfields} and \ref{groupentries}.
Furthermore, the quasi-relators are relations in these groups.
\end{thm}

\begin{rmk}
 It is proved in \cite{lipyan} that the quasi-relators generate all the relations
for these groups and that the groups $G_i$ are discrete co-compact subgroups of
${\rm PSL}(2,\mathbb{C})$.
%\nospacepsl.
\end{rmk}
\begin{rmk}
 $f_2$ and $w_2$ give an efficient way to solve the word problem in these groups.
\end{rmk}

\begin{table}

\begin{tabular}{|c|p{3in}|p{1in}|}

\hline
Region &Minimal Polynomial & Numerical Value \\
\hline
$X_0$ & $\tau^8 + 2\tau^6 + 6\tau^4 + 2\tau^2 + 1$ & $0.853230697-1.252448658i$ \\

\hline
$X_1$ & $ \tau^8 - 2\tau^7 + 5\tau^6 - 4\tau^5 + 7\tau^4 - 4\tau^3 + 5\tau^2 - 2\tau + 1 $ &
$0.904047196- 1.471654224i$ \\

\hline
$X_2$ & $\tau^4 - 2\tau^3 + 4\tau^2 - 2\tau + 1 $ &
$0.742934136 -1.529085514i$ \\

\hline
$X_3$ & $\tau^{24} - 8\tau^{23} + 35\tau^{22} - 107\tau^{21} + 261\tau^{20} - 538\tau^{19} +
972\tau^{18} - 1565\tau^{17} + 2282\tau^{16} - 3034\tau^{15} + 3706\tau^{14} - 4171\tau^{13} +
4339\tau^{12}- 4171\tau^{11} + 3706\tau^{10} - 3034\tau^9 + 2282\tau^8-1565\tau^7 + 972\tau^6 -
538\tau^5 + 261\tau^4 -107\tau^3 + 35\tau^2 - 8\tau + 1 $ &
$ 1.404292212 - 1.179267298i$ \\

\hline
$X_4$ & $\tau^6 - 3\tau^5 + 5\tau^4-4\tau^3 + 5\tau^2 - 3\tau + 1 $
 & $1.354619901-1.225125454i$ \\

\hline

   $X_{5}$ & $\tau^{12} + 2\tau^{10} + 7\tau^8 - 4\tau^6+ 7\tau^4 + 2\tau^2 + 1 $
& $0.868063287- 1.460023666i$ \\
\hline

   $X_{6}$ & $\tau^{12} - 2\tau^{10} + 7\tau^8 + 4\tau^6 + 7\tau^4 - 2\tau^2 + 1$
& $1.460023666- 0.868063287i$ \\

\hline
\end{tabular}
\vspace{.1in}
\caption{ Field containing $z$ for all regions}
\label{groupfields}
\end{table}

\begin{table}
\begin{tabular}{|c|p{2.5in}|c|p{1in}|}
\hline
Region & ${\rm tr_1}$ & ${\rm tr_2}$ & ${\rm tr_3}$ \\

\hline
$X_0$ & $-z-6z^3-2z^5-z^7$ & ${\rm tr_1}$ & $(-5z^2-2z^4-z^6)/2$ \\

\hline
$X_1$ & $2-4z+4z^2-7z^3+4z^4-5z^5+2z^6-z^7$ & ${\rm tr_1}$ & ${\rm tr_1}$ \\

\hline

$X_2$ & $2-3z+2z^2-z^3$ & ${\rm tr_1}$ & ${\rm tr_1}$\\

\hline
$X_3$ & $8-34z+107z^2-261z^3+538z^4-972z^5+1565z^6-2282z^7+3034z^8-3706z^9+4171z^{10}-4339z^{11}+
4171z^{12}-3706z^{13} +3034z^{14}-2282z^{15}+1565z^{16}-972z^{17}+ 538z^{18}-261z^{19}+107z^{20}-
35z^{21}+8z^{22}-z^{23}$ & ${\rm tr_1}$ & ${\rm tr_1}$ \\
\hline

$X_4$ & $3-4z+4z^2-5z^3+3z^4-z^5$ & ${\rm tr_1}$ & ${\rm tr_1}$ \\

\hline

   $X_{5}$ & $-z-7z^3+4z^5-7z^7-2z^9-z^{11}$ & ${\rm tr_1}$ & $(-6z^2+4z^4-7z^6-2z^8-z^{10})/2$ \\
\hline

   $X_{6}$ & $3z-7z^3-4z^5-7z^7+2z^9-z^{11}$ & ${\rm tr_1} $ & $(4-6z^2-4z^4-7z^6+2z^8-z^{10})/2 $ \\
\hline

\end{tabular}
\vspace{.1in}
\caption{Group coefficients as polynomials in $z$ in respective field}
\label{groupentries}
\end{table}

In this way we also obtain $\tr(f^2), \tr(w^2)$ and $\tr(f^2w^2)$ as roots of polynomials and
find a primitive element which generates the field that contains all the three traces. We have the following theorem. 
\begin{thm} The invariant trace fields for all the regions are as given in Table \ref{itf}.
\end{thm}

\begin{table}
\begin{tabular}{|c|p{3in}|p{1in}|}

\hline

Region & Minimal Polynomial & Numerical Value \\

\hline
$X_0$ & $\tau^2 + 3$ & $1.732050808i$ \\

\hline
$X_1$ & $ \tau^4 -2 \tau^3 + \tau^2 -2 \tau + 1$ & $-0.207106781 + 0.978318343i$ \\
\hline
$X_2$ & $\tau^2 + 1$ & $i$ \\
\hline
$X_3$ & $\tau^{12}+ 6\tau^{11} + 23\tau^{10} + 91\tau^9 + 257\tau^8 + 489\tau^7 + 823\tau^6 +
1054\tau^5 - 13\tau^4 - 2445\tau^3 - 3405\tau^2 - 1847\tau - 337$ &
$ 0.632778000- 3.019170376i $ \\

\hline
$X_4$ & $\tau^3 - \tau - 2$ & $-0.760689853 + 0.857873626i$ \\

\hline

    $X_{5}$ & $\tau^3-\tau^2+\tau+1$ & $0.771844506 + 1.11514250i$ \\
\hline

  $X_{6}$ & $\tau^3-\tau^2+\tau+1$ & $0.771844506 - 1.11514250i$ \\
\hline

\end{tabular}
\vspace{.1in}
\caption{Invariant trace fields for all the regions}
\label{itf}
\end{table}

\begin{rmk}
The invariant trace field descriptions in Table \ref{itf} for $X_i$ for $i\neq 3$ are the
canonical field description given by Snap (see \cite{cghn00}).
\end{rmk}

 \section{The manifolds for the regions $X_1, X_2,X_5$ and $X_6$}

 In this section we find manifolds from the Hodgson \& Weeks census of closed
hyperbolic 3-manifolds whose fundamental groups are isomorphic to the groups $G_i$ for
$i=1,2,5,6$. This census is included in Jeff Weeks' program SnapPea \cite{snappea}
and is refered to as SnapPea's census of closed hyperbolic 3-manifolds. These manifolds
are described as Dehn surgeries on cusped hyperbolic 3-manifolds from SnapPea's census
of cusped manifolds  \cite{chw99}, \cite{hw89}. We use the invariant
trace fields and volume estimates for the
regions given in \cite{jr01} and \cite{lipyan} to search through the roughly 11,000
manifolds in the closed census. The package
Snap \cite{snap} includes a text file called closed.fields which lists the
invariant trace fields for the manifolds in the closed census.
Using this file to compare the invariant trace fields we narrowed our search
to less than 50 manifolds for each regions. Then using the homology,
volume estimates and length of shortest geodesic we further narrowed the
search to less than 5 manifolds. Table \ref{mfddata} gives the
approximate volume ($V$) as given in \cite{jr01}, first homology ($H_1$),
approximate length of shortest geodesic which is the values of the
parameter $L$ ($l_{min}$) and the manifold description as given in SnapPea.

\begin{table}
\begin{tabular}{|l|c|c|c|p{2in}|}
\hline
Region & $V$ & $H_1$ & $l_{min}$ & Manifolds \\
\hline

$X_1$ & 4.11696874 & $\ZZ_7 \oplus \ZZ_7$ & 1.0930 &
$v2678(2,1),v2796(1,2)$ \\

\hline
 $X_2$ & 3.66386238 & $\ZZ_4 \oplus \ZZ_{12}$ & 1.061 &
$s778(-3,1), v2018(2,1)$ \\

\hline

 $X_4$ & 7.517689 & $\ZZ_4 \oplus \ZZ_{12}$ & 1.2046 & NA \\

\hline

 $X_5$ or $X_6$ & 3.17729328 & $\ZZ_4 \oplus \ZZ_4$ & 1.0595 &
$s479(-3,1), s480(-3,1)$, $s645(1,2), s781(-1,2)$, $v2018(-2,1)$ \\

\hline

\end{tabular}
\vspace{.1in}
\caption{Data for regions $X_1,X_2,X_4,X_5$ and $X_6$}
\label{mfddata}
\end{table}

The above manifolds include the manifolds mentioned in \cite{gmt}
for the regions $X_1$, $X_2$ and $X_5$. All the SnapPea manifolds
associated to a region in Table \ref{mfddata} are isometric. It is
shown in \cite{jr01} that the manifolds associated to the regions $X_5$
and $X_6$ are isometric with an orientation reversing isomtery. The
manifold associated to $X_4$ is discussed in the next section and
that for $X_3$ is discussed in \cite{lipyan}.

\begin{table}
\begin{center}
\begin{tabular}{|l|p{3in}|}
\hline Manifold & $\pi_1$ Relators \\
\hline $v2678(2,1)$

& $q_1=a^2b^2aba^{-1}ba^{-1}b^{-1}a^{-1}ba^{-1}bab^2$\\
& $q_2=ab^{-1}ab^{-1}a^{-1}b^{-1}ab^{-1}aba^2b^2a^2b$\\

\hline $s778(-3,1)$ &

$q_1=ab^{-1}aba^2b^2ab^2a^2bab^{-1}$ \\
& $q_2=ab^2a^2ba^2b^2aba^{-1}ba^{-1}b$\\

\hline $s479(-3,1)$ &

$q_1=aba^2b^2a^2bab^{-2}a^{-2}b^{-2}$ \\
& $q_2=a^2b^2ab^2a^2bab^{-1}ab^{-1}ab$\\
\hline

\end{tabular}
\end{center}
\vspace{.1in}
\caption{Relators for manifolds}
\label{relatorsmfd}
\end{table}

\begin{table}
\begin{center}
\begin{tabular}{|l|l|l|l|}
\hline Region & Manifolds & Isomorphism & Inverse\\

\hline $X_1$ &$v2678(2,1)$ & $f\longrightarrow a^{-1}$, $w \longrightarrow b$ &

$a \longrightarrow f^{-1} $, $ b \longrightarrow w$ \\
\hline
$X_2$ &$s778(-3,1)$&

$f \longrightarrow a$, $w \longrightarrow b^{-1}$ & $a \longrightarrow f$,

$b \longrightarrow w^{-1}$ \\
\hline

$X_5$ & $s479(-3,1)$ &

$f \longrightarrow ab$,

$w \longrightarrow b$ & $a \longrightarrow fw^{-1}$, $b \longrightarrow w$ \\
\hline
\end{tabular}

\end{center}
\vspace{.1in}
\caption{Isomorphisms}
\label{isomorphisms}
\end{table}

The fundamental groups of the above manifolds have two generators and two
relators. Table \ref{relators} and \ref{relatorsmfd} gives the relators
for the marked groups and the relators for the fundamental
groups of the corresponding manifolds respectively. One can verify the isomorphisms
between the groups given in Table \ref{isomorphisms}. We have the following theorem.
\begin{thm}
The manifolds $v2678(2,1),s778(-3,1)$ and $s479(-3,1)$ in SnapPea's census of closed
manifolds are exceptional manifolds associated to the regions $X_1, X_2$ and $X_5$
respectively.
\end{thm}

 \section{The manifold associated to the region $X_4$}

In this section we give a description of the manifold associated to the
region $X_4$ as a
double cover of an orbifold commensurable to the manifold $m369(-1,3)$ in
SnapPea's census of closed manifolds.

In Section $4$ using the approximate volumes and other data given in \cite{jr01} and \cite{lipyan}
we found manifolds from the SnapPea's census of closed manifolds
with fundamental groups isomorphic to the groups for the regions $X_1, X_2, X_5$ and $X_6$.
The regions $X_3$ and $X_4$ could not be handled because of their large
volumes. However for
the region $X_4$ a list of manifolds was found in the closed census having
approximately half the volume of
$X_4$ and the same commensurability invariants.

% These manifolds are: $s297(1,3)$,
% $s298(5,1)$, $s594(1,2)$, $m307(-5,1)$, $m305(-5,1)$, $m369(-1,3)$,
% $m371(1,3)$,
% $m290(-1,4)$, $m390(3,1)$, $m293(-2,3)$, $m303(1,3)$, $s594(2,1)$,
% $s480(3,1)$,
% $s595(1,2)$, $s235(-4,3)$, and $s287(-3,1)$.

In hope of obtaining the manifold for $X_4$ as a double cover of one
of these manifolds we compared index two subgroups of the fundamental
groups of each of these manifolds to $G_4$, the marked group for
$X_4$. Most of the subgroups were eliminated on basis of homology,
however one index two subgroup of the census manifold $m369(-1,3)$ had
the correct homology, and same lengths for its elements as for
$X_4$. Using the program \textit{testisom} \cite{testisom} it was
checked that this subgroup was not isomorphic to $G_4$.

% However using the similarity of the
% geometric information for this subgroup and $G_4$ we obtained an orbifold
% quotient of this subgroup which had the $G_4$ as an index two subgroup. The
% computations and presentations were obtained using {\it GAP} and {\it
% magnus} and the
% isomorphisms were checked using {\it testisom}.\\
% \newline

\begin{thm}
The manifold $N_4$ associated to the region $X_4$ is commensurable with the manifold
$m369(-1,3)$ in SnapPea's census of closed manifolds. This manifold is
obtained as a double cover of an orbifold which is double covered by a double cover of $m369(-1,3)$. \end{thm}

\pf Let $M=m369(-1,3)$. We will construct the following diagram of $2:1$ covers: $$ \begin{array}{cccc}
N & & N_4 \\
\downarrow & \searrow & \downarrow \\
M & & O
\end{array}
$$
We obtain a presentation of $\pi_1(M)$ from SnapPea.
$$ \pi_1(M)= \langle a,\ b,\ c|\ ab^{-1}a^{-1}c^2bc,\ abcb^3a^{-1}c^{-1},\ acbc^{-1}b^{-1}cbacb \rangle$$

Let $\phi:\pi_1(M)\rightarrow \mathbb{Z}_2$ be defined by
$\phi(a)=1,\phi(b)=\phi(c)=0$. Then $\phi$ is a homomorphism and ker($\phi$)
is an index two subgroup of $\pi_1(M)$ generated by $b$ and $c$. Let $N$
denote the double cover of $M$ corresponding to this subgroup so that
$\pi_1(N)= {\rm ker}(\phi)$. A presentation of $\pi_1(N)$ is
$$ \pi_1(N)= \langle b,\ c |\ r_1,\ r_2 \rangle$$
where
\begin{eqnarray*}
r_1&=& bcb^3cbc^{-1}b^{-1}cbc^{-1}b^{-1}cbc^2(bc^3)^2bc^2bcb^{-1}c^{-1}bcb^{-1}c^{-1}, \\
r_2&=& cbc^{-1}b^{-1}cbc^{-1}b^{-1}cbc^2(bc^3)^2(bc^2bc^3bc^3bc^3)^2bcb^{-1}c^{-1}bcb^{-1}c^{-1}b.
\end{eqnarray*}

Let $\psi: \pi_1(N)\rightarrow \pi_1(N)$ be defined by $\psi(c)= c^{-1}$ and
$\psi(b)= c^3b$. Then $\psi$ is an automorphism of $\pi_1(N)$ of order
two. Extending the group $\pi_1(N)$ by this automorphism we obtain a group
$H$ whose presentation is
$$H=\langle b,\ c,\ t |\ r_1,\ r_2,\ tct^{-1}c,\ tbt^{-1}b^{-1}c^{-3},\ t^2\rangle$$

$\pi_1(N)$ is a subgroup of $H$ of index two and the quotient of
$\mathbb{H}^3$ by $H$ is an orbifold $O$ (due to the torsion element $t$)
which is double covered by $N$. Let $\mu: H \rightarrow \mathbb{Z}_2$ be
defined by $\mu(c)=0,\ \mu(b)=\mu(t)=1$. Then $\mu$ is a homomorphism and
ker($\mu$) is an index $2$ subgroup of $H$ generated by elements $c$ and $b*t$.
Let $x=c$ and $y=b*t$. Then a presentation of ker($\mu$)$=G$ is
$$G= \mathrm{ker}(\mu)=\langle x,\ y |\ s_1,\ s_2,\ s_3 \rangle$$
where
\begin{eqnarray*}
s_1&=& (yx^{-1}y^{-1}x^{-1})^2yx^2y^2x^3yxyxy^{-1}xyxy^{-1}xyxyx^3y^2x^2, \\
s_2&=& (yxy^{-1}x)^2yxyx^3(y^2x^2y^2x^3)^2yxyxy^{-1}xyxy^{-1}xyx, \\
s_3&=& y^{-1}x^{-3}(y^{-1}x^{-1})^2yx^{-1}y^{-1}x^{-1}yx^2y^2x^3yxyx^3y^2x^2(yx^{-1}y^{-1}x^{-1})^2.
\end{eqnarray*}

 The presentation for the marked group $G_4$ as given in \cite{gmt} is
$$G_4=\langle f,\ w |\ r_1(X_4),\ r_2(X_4) \rangle$$
where
\begin{eqnarray*}
r_1(X_4)&=& f^{-2}wfwf^{-1}(wfw^{-1}f)^2wf^{-1}wfwf^{-2}w^2(f^{-1}w^{-1}f^{-1}w)^2w, \\
r_2(X_4)&=&f^{-1}(f^{-1}wfw)^2f^{-2}w^2f^{-1}w^{-1}f^{-1}w(f^{-1}w^{-1}fw^{-1})^2 \\
& & \times \ f^{-1}wf^{-1}w^{-1}f^{-1}w^2.
\end{eqnarray*}

One easily verifies that the map $\nu:G_4 \rightarrow G$ given by
$\nu(x)=f$ and $\nu(y)=f^{-1}w^{-1} $ is an isomorphism. The inverse
of $\nu$ is given by $\nu^{-1}(f)=y $ and
$\nu^{-1}(w)=y^{-1}x^{-1}$. Lipyanskiy \cite{lipyan} constructed a
Dirichlet domain for all the regions whose groups are isomorphic to
the marked groups. It follows that $G_4$ is torsion free and hence $G$
is a torsion free subgroup of $H$ of index $2$. Hence it gives the
manifold $N_4$ which double covers the orbifold $O$.\done
\begin{rmk}
The symmetries of the configuration of lines in \h consisting of the
axis of $f$, $w$, their translations and the orthocurves between them
led us to study the above mentioned subgroups and automorphisms.
%We began looking more closely at
$N$ has a geodesic of the same length and the same translation
length as $N_4$ but it is not the shortest geodesic
in $N$.
\end{rmk}

 \section{Uniqueness}

 In this section we address the issue of uniqueness of the solutions
in the  given region.
In Chapter $3$ of \cite{gmt} the authors showed
the existence and uniqueness of solution for the region $X_0$ by
using a geometric argument to show $R'=1$ and then using the symmetry
of the region $X_0$ to reduce the number of variables and obtain a one
variable equation which has only one solution in the region $X_0$.
Using Groebner basis we show that there is a unique point in every
region $X_i$ for which the quasi-relators equal the identity.

%We show that the solutions we found are the only solutions in  the
%given region for all the regions.

Let $I$ be the ideal generated by the equations formed by the entries of
the quasi-relators of a region subtracted from the identity matrix.
We compute a Groebner basis for $I$ and verify that there is a
unique solution to equations in the Groebner basis in that region.
% we repeat this for each of the variables by  using different
%orderings.
For computational convenience we split the  relations in
half as it reduces the degree of the polynomials generating  the
ideal.
%We use the matrix representatives given in Section 3,  Equation
%\ref{fandwnew} for $f$ and $w$ instead of the ones given in  Section 2
%since this also reduces the degree of the polynomials.
Let $p={\rm tr_1}=\tr(f)$, $q={\rm tr_2}=\tr(w)$ and $r={\rm tr_3}=\tr(f^{-1}w)$ as in
Section 3. Using Equation \ref{fandw} in Section 2, $L',\ D'$ and $R'$ can
be expressed in terms of $p,\ q$ and $r$ as follows:
\begin{eqnarray}
L'&=& \Big( \frac{p\pm \sqrt{p^2-4}}{2}\Big) ^2 \\
D'&=&\Big(\frac{2q\sqrt{R'}\pm \sqrt{4q^2R'-4(1+R')^2}}{2(1+R')}\Big)^2\\
R'&=& \frac{qL'-r\sqrt{L'}}{r\sqrt{L'}-q}
\end{eqnarray}
\noindent
Using Equation \ref{fandwnew} in Section 3 we can write the entries of
conjugates of $f$ and $w$ in terms of $p,\ q,\ r$ and $z$ where $qz-z^2=1$.
The equations for quasi-relators using Equation \ref{fandwnew} are simpler
for computing Groebner basis.

For example, for the region $X_0$, using the ordering $z,\ r,\ q,\ p$
on the variables, the last entry of the Groebner basis is
$(p-1)(p+1)(p^4-2p^2+4)$.
Using Equation 6 it can be easily checked that only one root of the above equation
gives the value of $L'$ lying in the region $X_0$.
 Similarly the last entries of the Groebner bases in
orders $z,\ r,\ p,\ q$ and  $z,\ p,\ q,\ r$ are
$(q-2)(q+2)(q^4-2q^2+4)$ and $(r+1)(r^2-r+1)$ respectively. Using Equations
7 and 8 it can be easily checked that only one root of $q^4-2q^2+4$ and $r^2-r+1$
give the value of $D'$ and $R'$ respectively lying in the region $X_0$.
This shows that there is a unique solution for the quasi-relators
in the region $X_0$.

Similarly for the regions $X_2,\ X_4,\ X_5$ and $X_6$ the
last entry of
the Groebner basis is a polynomial in either $p,\ q$ or $r$
depending on the ordering of the variables. Using Equations 6, 7 and 8
we check that there is a unique solution in the respective region.

For the regions $X_1$ and $X_3$ we obtain a multivariable polynomial in
$p$, $q$ and $r$ as a factor of the last entry of the Groebner basis along
with a single variable polynomial. We eliminate this factor using the
Mean Value Theorem.

For example, for the region $X_3$ the last entry of the Groebner basis
with the ordering $z>r>q>p$ on the variables factors as:

$(p^3 + p^2 - 2p - 1)
     (p^{10} - 7p^9 + 15p^8 + 4p^7 - 49p^6 + 11p^5 + 88p^4
  + 87p^3 - 501p^2 + 543p - 193)
     (p^{10} + 5p^9 + 6p^8 - 6p^7 - 10p^6 + 12p^5 + 13p^4
 - 11p^3 - 6p^2 + 4p + 1)
     (p^{12} - 8p^{11} + 23p^{10} - 19p^9 - 35p^8 + 73p^7 - 3p^6
 - 72p^5 + 25p^4 + 29p^3 - 11p^2 - 3p + 1)
     (rq + rp - r + qp - q - p + 1)$

It can be checked that only one root of the polynomial in $p$ above gives
the value of $L'$ lying in the region $X_3$.
%is obtained
%using Equation 6 from a root of the polynomial $p^{12} - 8p^{11} + 23p^{10} - 19p^9 -
%35p^8 + 73p^7 - 3p^6 - 72p^5 + 25p^4 + 29p^3 - 11p^2 - 3p + 1$. The roots of the
%other polynomials in $p$ give values of $L'$ which do not lie in the region $X_3$.
We will show that the polynomial $rq + rp - r + qp - q - p + 1$ has no root
in the region $X_3$.

From Table \ref{rangeforx3}, we see that
\begin{eqnarray*}
L'&=& 0.581385 - 3.312055\, i, \\
D' &=& 1.15663 - 2.756005\, i, \\
R' &=& 1.40437 - 1.179435\, i
\end{eqnarray*}
is the midpoint of region $X_3$.
Then using Equation \ref{fandw} the $p,\ q,\ r$ values corresponding to this
point are
\begin{eqnarray*}
p_0&=&1.8219- 0.828571\, i,\\
q_0 &= & 1.82191 - 0.828633\, i,\\
r_0 &=& 1.82192 - 0.828537\, i
\end{eqnarray*}
If the polynomial $f(p,q,r)=rq + rp - r + qp - q - p + 1$  has a root,
say $(p_1,q_1,r_1)$, in the
region $X_3$ then by the Mean Value Theorem for some point $(p,q,r)$ in $X_3$
we obtain
\begin{eqnarray*}
|f(p_0,q_0,r_0)|&=&|\nabla f(p,q,r)\cdot (p_1-p_0,q_1-q_0,r_1-r_0)|\\
& \leq& ||\nabla f(p,q,r)|| \, ||(p_1-p_0,q_1-q_0,r_1-r_0)||
\end{eqnarray*}
From the parameter ranges of region $X_3$ from Table \ref{rangeforx3}
we know that $||(p_1-p_0,q_1-q_0,r_1-r_0)||< 0.002$.
Hence $||\nabla f(p,q,r)|| \geq |f(p_0,q_0,r_0)|/0.002$ at some point in the region
$X_3$. It can be checked that $|f(p_0,q_0,r_0)|/0.002 \sim 3000$ and that
\begin{eqnarray*}
||\nabla(f)|| &\leq& \sqrt{(|p|+|q|+1)^2+(|q|+|r|+1)^2+(|r|+|p|+1)^2} \\
% &\leq&\sqrt{(|p_0|+|q_0|+2h+1)^2+(|q_0|+|r_0|+2h+1)^2+(|r_0|+|p_0|+2h+1)^2}\\
 &\leq&\sqrt{(|p_0|+|q_0|+2)^2+(|q_0|+|r_0|+2)^2+(|r_0|+|p_0|+2)^2}\\
 &<&11
\end{eqnarray*}
\noindent
in the region $X_3$ and hence $f(p,q,r)$ does not have a root in $X_3$.

We can similarly check for the other variables by changing the order of the
variables. The
last entry of the Groebner basis with the ordering $z>r>p>q$ on the variables factors as

$    (q - 2)
    (q^{10} - 7q^9 + 15q^8 + 4q^7 - 49q^6 + 11q^5 + 88q^4
+ 87q^3 - 501q^2 + 543q - 193)
    (q^{10} + 5q^9 + 6q^8 - 6q^7 - 10q^6 + 12q^5 + 13q^4
- 11q^3 - 6q^2 + 4q + 1)
    (q^{12} - 8q^{11} + 23q^{10} - 19q^9 - 35q^8 + 73q^7 - 3q^6
- 72q^5 + 25q^4 + 29q^3 - 11q^2 - 3q + 1)
    (rp + rq - r + pq - p - q + 1)$

\noindent
and the last entry of the Groebner basis with the ordering $z>p>q>r$ on the variables factors as

$    (r^3 + r^2 - 2r - 1)
    (r^{10} - 7r^9 + 15r^8 + 4r^7 - 49r^6 + 11r^5 + 88r^4
+ 87r^3 - 501r^2 + 543r - 193)
    (r^{10} + 5r^9 + 6r^8 - 6r^7 - 10r^6 + 12r^5 + 13r^4
- 11r^3 - 6r^2 + 4r + 1)
    (r^{12} - 8r^{11} + 23r^{10} - 19r^9 - 35r^8 + 73r^7 - 3r^6
- 72r^5 + 25r^4 + 29r^3 - 11r^2 - 3r + 1)
    (pq + pr - p + qr - q - r + 1)$

Using Equations 6, 7 and 8 above it can be checked that there is only one
root of the above polynomials in $q$ and $r$ which give values for
$D'$ and $R'$ respectively lying in the region $X_3$.
 %are obtained from a root
%of the polynomials $q^{12} - 8q^{11} + 23q^{10} - 19q^9 - 35q^8 +
%73q^7 - 3q^6  - 72q^5 + 25q^4 + 29q^3 - 11q^2 - 3q + 1$ and $r^{12} -
%8r^{11} + 23r^{10} - 19r^9 - 35r^8 + 73r^7 - 3r^6  - 72r^5 + 25r^4 +
%29r^3 - 11r^2 - 3r + 1$.  The roots of the other polynomials in $q$
%and $r$ give values of $D'$ and $R'$ which do not lie in the region
%$X_3$.
The polynomial in $p, \ q$ and $r$ is the same as above. Hence
the quasi-relators for the region $X_3$ have a unique solution in the
region $X_3$. The uniqueness of solutions is proved similarly for
the region $X_1$.
\begin{prop}\label{unique}
Let $f$ and $w$ be as in Equation 1 and let $r_1(X_i), \ r_2(X_i)$
be quasi-relators for the region $X_i$. Then there is a unique triple $(L',D',R')$
in the region $X_i$ for which the quasi-relators equal the identity matrix.
\end{prop}
\pf It follows from Theorem 3 that there is a triple $(L',D',R')$ in
the region $X_i$ for which the quasi-relators equal the identity
matrix. The uniqueness follows from the Groebner basis computation for
every region. \done

%Let the hyperbolic distance between two points $x$ and $y$ of \h be
%denoted by $\rho(x,y)$. For an isomtery $f$ of \h let
%\textit{Relength(f)}$=\mathrm{inf}\{\rho(x,f(x))|x\in \h \}$.
%Let
%$N_0=$Vol3, $N_1=v2678(2,1),\ N_2=s778(-3,1)$, $N_3$ be the
%exceptional manifold associated to $X_3$ found in \cite{lipyan}, $N_4$
%be exceptional manifold associated to $X_4$ found in Section 5 and
%$N_5=s479(-3,1)$. Then Theorem 1 can be stated as

%\begin{thm}Let $N$ be an exceptional manifold. Then $N$ is covered by $N_i$ for  some $i=0,1,2,3,4,5$.
%\end{thm}
We now give the proof of our main Theorem. \\
\ \\
{\bf Proof of Theorem 1:}
Let $N$ be an exceptional manifold and let $\delta$ be the shortest
geodesic in $N$. Let $f \in \pi_1(N)$ be a primitive hyperbolic
isometry whose fixed axis $\delta_0 \in \mathbb{H}^3$ projects to
$\delta$ and $w\in \pi_1(N)$ be a hyperbolic isometry which takes
$\delta_0$ to its nearest translate. Let $G$ be the subgroup
of $\pi_1(N)$ generated by $f$ and $w$. Then the manifold
$N'=\h/G$ is exceptional and $\delta_0$ projects to the
shortest geodesic in $N'$.

It follows from Corollary 1.29 of \cite{gmt} that the
$(L',D',R')$ parameter for $G$ lies in the region $X_i$ for
some $i=0,1, \ldots ,6$. By definition of the quasi-relators
\cite{gmt},
$Relength(r_1),\ Relength(r_2) < Relength(f)$. Since $f$ is
the shortest element in $G$, $r_1$ and $r_2$ equal the identity
in $G$ i.e. they are relations in $G$.  It is proved in
\cite{lipyan} that the quasi-relators generate all the relations
in the groups $G_i= \langle f, w| r_1(X_i), r_2(X_i)\rangle$
for $i=0, \ldots, 6$. Hence $G=G_i$
 and $\pi_1(N)$ contains the marked group $G_i$ for some $i=0,1,
\ldots 6$.

It follows from Proposition 1 that the quasi-relators for a given
region equal the identity at a unique point in that region. Hence
$N$ is covered by $N_i$ for some $i=0,1,2,3,4,5$ where
$N_i$ are the manifolds described in the Introduction. \done
\section{Conclusions}
%%  Let $T$ consists of those parameters corresponding to the groups
%%  $\{G,\ f,\ w\}$ such that $Relength(f)$ is the shortest element of
%%  $G$ and the distance between the axis of $f$ and its nearest
%%  translate is less than $\ln(3)$. Let $S=exp(T)$. In \cite{gmt} the
%%  authors conjecture:
%% \begin{conj}Each sub-box $X_i$, $0\leq i \leq 6$ contains a unique element $s_i$ of S.
%%  Further if $\{G_i,f_i,w_i\}$ is the marked group associated to $s_i$ then $N_i=\h/G_i$
%% is a closed hyperbolic 3-manifold with the following properties
%% \begin{enumerate}
%% \item[(i)] $N_i$ has fundamental group $\langle f,\ w | r_1(X_i),\
%% r_2(X_i) \rangle $ where $r_1$ and $r_2$ are the quasi-relators
%% associated to the box $X_i$.
%% \item[(ii)]$N_i$ has a Heegaard genus 2 splitting realizing the above
%% group presentation.
%% \item[(iii)]$N_i$ nontrivially covers no manifold.
%% \item[(iv)]$N_6$ is isometric to $N_5$.
%% \item[(v)]If $(L_i,D_i,R_i)$ is the parameter in $T$ corresponding to
%% $s_i$, then $L_i,\ D_i,$ $R_i$ are related as follows: \\ For $X_0, \
%% X_5,\ X_6$: $L=D$, $R=0$.\\ For $X_1,\ X_2,\ X_3,\ X_4$: $R=L/2$.
%% \end{enumerate}
%
%% \end{conj}
%% In Chapter $3$ of \cite{gmt} the authors settle this conjecture for
%% the region $X_0$. Jones and Reid \cite{jr01} prove part (i), (iii) for
%% region $X_0$ and (iv).
%
The first part of Conjecture \ref{gmtconj} follows from Proposition \ref{unique}.
The results in Section $4$ and $5$ shows part (i) of Conjecture \ref{gmtconj}
 for regions $X_i$, $i= 1,2,4,5,6$ and the exact arithmetic from Section $3$
show part (v) of Conjecture \ref{gmtconj} for all the regions.
The question about the
uniqueness of the manifolds remains open for all regions except
$X_0$. It is reasonable to make the following  conjecture:
\begin{conj}
The manifolds $N_i$ for $i=1,2,3,4,5,6$ do not nontrivially cover any manifold.
\end{conj}
Alan Reid proves the conjecture for $N_1$ and $N_5$ in the following appendix.

\section{Appendix: The manifolds $N_1$ and $N_5$}

\centerline{\bf Alan W. Reid}

In this appendix we prove the following theorem.

\medskip

\begin{thmar}
The manifolds $N_1$ and $N_5$ do not properly
cover any closed hyperbolic 3-manifold
\end{thmar}

\pf
We give the proof in the case of $N_1$, the case
of $N_5$ is similar. Both arguments follow that given for $\Vol3$
in \cite{jr01}. We refer the reader to \cite{mr03} for details about arithmetic Kleinian
groups and quaternion algebras.

Thus, suppose that $N_1 = {\Bbb H}^3/\Gamma_1$
non-trivially covered
a closed hyperbolic 3-manifold $N= {\Bbb H}^3/\Gamma$ say, with covering
degree $d$. Using the identification of $N_1$ as $v2678(2,1)$
given in the paper,
it follows that the volume of $N_1$ is approximately
$4.116968736384613\ldots$ and $H_1(N_1;{\Bbb Z})= \ZZ/7\ZZ \oplus \ZZ/7\ZZ$ is finite of odd order.
Note that since $H_1(N_1;{\Bbb Z})$ is finite
the closed hyperbolic 3-manifold $N$
is orientable.
Since the volume of the smallest arithmetic manifold is
approximately $0.94$ \cite{cfjr01}, it follows that $d\leq 4$.

Using Snap (or from computations in the paper) the Kleinian group
$\Gamma_1$, and hence $\Gamma$, is arithmetic with invariant
trace-field $k$ say. This has degree 4, has discriminant $-448$, and the
invariant quaternion algebra $B/k$ unramified at all finite places. We
remark that there is a unique such field.

Since $|H_1(N_1;{\Bbb Z})|$ is odd, $\Gamma_1$ is
derived from a  quaternion algebra. Furthermore, $B$ has type number
1, and so $\Gamma_1$ is conjugate into the group of elements of norm 1
of a maximal order $\mathcal O$ of $B$. Now the image of the elements
of norm 1 of $\mathcal O$ in
$\PSL(2,{\Bbb C})$ can be shown (see \cite{mr89}) to coincide with the orientation
preserving subgroup of index 2 in the Coxeter simplex group
$T[2,3,3;2,3,4]$. The notation for the Coxeter group is that of \cite{mr89}.

Denote this group by $C$. The minimal index of a torsion free subgroup
of $C$ is at least $24$, since by inspection of the Coxeter
diagram, $C$ contains a subgroup isomorphic to $S_4$.  Therefore the volume
calculations of \cite{mr89} show that $\Gamma_1$ is a minimal index
torsion-free subgroup of this group.

Now the analysis in \S 4 of \cite{jr01}
shows that the possible maximal groups in the
commensurabilty class of $\Gamma_1$ that contain $\Gamma$ are either
the group $\Gamma_{\mathcal O}$ (in the notation \cite{jr01}) where $\mathcal O$ is
the maximal order above, or $\Gamma_{\{{\mathcal P}_7\},{\mathcal O}}$ (in the
notation of \cite{jr01}).

Suppose first that $\Gamma < \Gamma_{\mathcal O}$. By the remarks above
$\Gamma$ is not a subgroup of $C$. Now \cite{mr89} shows that $\Gamma_{\mathcal
  O}$ contains $C$ as a subgroup of index 2.  It follows that $\Gamma$
must contain $\Gamma_1$ as a subgroup of index $2$, and that $\Gamma$ is a
torsion-free subgroup of index $24$ in $\Gamma_{\mathcal O}$.
However we claim that this is impossible.

Firstly we can obtain a presentation of the group
$\Gamma_{\mathcal O}$ using the geometric description of $C$ above; namely
the group $\Gamma_{\mathcal O}$ is obtained from $C$ by adjoining an
orientation-preserving involution $t$ that is visible in the Coxeter
diagram.  On checking the action of this involution, one gets that a
presentation for $\Gamma_{\mathcal O}$ is given by: \\
\ \\
$<t,a,b,c\ |\ t^2,\ a^2,\ b^3,\ c^3,\ (b*c)^2,\ (c*a)^3,\ (a*b)^4,\
                           t*a*t^{-1}*c*b,\ t*b*t^{-1}*a*c^{-1},\ t*c*t^{-1}*c>.$
\\
\ \\
Now a check with Magma (for example) shows that there are $24$
subgroups of index $24$ and all but two are easily seen to have
elements of finite order by inspection of presentations. The remaining
two have abelianizations ${\Bbb Z}/22{\Bbb Z}$. The index 2 subgroups in
these groups all have ${\Bbb Z}/11{\Bbb Z}$ in their abelianizations by
a standard cohomology of groups argument (or further checking with Magma). In
particular these index 2 subgroups cannot coincide with $\Gamma_1$,
which completes the analysis in this case.

For the second case the covolume of $G=\Gamma_{\{{\mathcal P}_7\},{\mathcal
O}}$ can be computed (using the formula in \S 2 of \cite{jr01}) to be 12
times that of $\Gamma_1$. An alternative, equivalent description of
this maximal group, is as the normalizer of an Eichler order $\mathcal E$ of level
${\mathcal P}_7$ in $\mathcal O$ (see \cite{mr03} Chapter 11 for example).

The results of \cite{cf00} (see in particular Theorems 3.3 and 3.6)
show that $G$ contains elements of orders $2$ and $3$, and so the
minimal index of a torsion-free subgroup in $G$ is at least 6.

Now $\Gamma_1 \subset G\cap C$. Furthermore, if we denote the image of
the group of elements of norm 1 in $\mathcal E$ in $\PSL(2,{\Bbb C})$ by
$\Gamma_{\mathcal E}^1$, then since the level is ${\mathcal P}_7$ it follows
that the index $[C:\Gamma_{\mathcal E}^1]$ is $8$ (see \cite{mr03} Chapters 6 and 11).
It is not difficult to see that
$G\cap C = \Gamma_{\mathcal E}^1$. One inclusion is clear, and the other
follows since, from above $[C:\Gamma_{\mathcal E}^1]=8$ so that the only
possible indices for $[C:G\cap C]$ are $2$ or $4$ (it cannot be 1 since
$G$ is a different maximal group from $\Gamma_{\mathcal O}$
above). Now $C$ is perfect, and so has no solvable quotients. Hence this rules out
$C$ from having index 2 or 4 subgroups.

%However $C$ is perfect, and so has no solvable quotients, and
%index 2 or 4 always provides such quotients.

We deduce from the above that $\Gamma_1$ is a subgroup of
$\Gamma_{\mathcal E}^1$.  Using the presentation of $C$, and on checking
with Magma for instance, we find that there are 5 subgroups of index
8, and some further low index computations on these subgroups using
Magma shows that only one can contain $\Gamma_1$. This subgroup (denoted
$H$ in what follows) is generated by two elements of order $3$ ( b and
ca in the generators above).

As in the first case we can use the geometry associated to $H$ to
construct a presentation for the group $G$. $H$ is 2-generator with
both generators of order $3$, so we can adjoin involutions $s$ and $t$
so that a presentation for $G$ is:
\\
\ \\
$<x,y,s,t\ |\ s^2,\ t^2,(s*t)^2,\ s*a*s*y^{-1},s*y*s^{-1}*x^{-1},t*x*t*x,
\ t*y*t^{-1}*y,\ x^3,\ y^3,\
x*y^{-1}*x^{-1}*y*x*y*x^{-1}*y^{-1}*x*y*x*y^{-1}*x^{-1}*y^{-1}*
x*y*x^{-1}*y^{-1}*x^{-1}*y*x*y^{-1}*x^{-1}*y^{-1}*x*y>.$
\\
\ \\
From our remarks above we need only check for torsion-free subgroups
of index $6$. However, an easy inspection using Magma shows that all
the subgroups of index $6$ (there are 4) have elements of finite
order. Hence this completes the proof. \done

\noindent
\textbf{Acknowledgments:} This work was supported in part
by grants from the NSF and the Texas Advanced Research Program.

\bibliography{ex_regions.bib}
\bibliographystyle{amsplain}

\end{document}